\documentclass[11pt]{amsart}
\usepackage{amssymb,latexsym}
\input epsf
\usepackage{graphicx}
\newtheorem{theorem}{Theorem}[section]
\newtheorem{proposition}[theorem]{Proposition}

\newtheorem{remark}[theorem]{Remark}
\newtheorem{lemma}[theorem]{Lemma}

\newtheorem{conjecture}[theorem]{Conjecture}

\newcommand{\cC}{{\mathcal C}}
\newcommand{\eE}{{\mathcal E}}

\newcommand{\gG}{{\mathcal G}}
\newcommand{\nN}{{\mathcal N}}

\newcommand{\KK}{{\mathbb K}}
\newcommand{\RR}{{\mathbb R}}

\newcommand{\sm}{{\smallsetminus}}
\begin{document}
\title[On the connectivity of polytope skeleta]
{On the graph-connectivity of skeleta of convex polytopes}

\author{Christos~A.~Athanasiadis}

\address{Department of Mathematics
(Division of Algebra-Geometry)\\
University of Athens\\
Panepistimioupolis\\
15784 Athens, Greece}
\email{caath@math.uoa.gr}

\date{January 10, 2008.}
\thanks{ 
Supported by the 70/4/8755 ELKE Research Fund of the University of
Athens.}
\begin{abstract}
Given a $d$-dimensional convex polytope $P$ and nonnegative
integer $k$ not exceeding $d-1$, let $\gG_k (P)$ denote the simple
graph on the node set of $k$-dimensional faces of $P$ in which two
such faces are adjacent if there exists a $(k+1)$-dimensional face
of $P$ which contains them both. The graph $\gG_k (P)$ is
isomorphic to the dual graph of the $(d-k)$-dimensional skeleton
of the normal fan of $P$. For fixed values of $k$ and $d$, the
largest integer $m$ such that $\gG_k (P)$ is $m$-vertex-connected
for all $d$-dimensional polytopes $P$ is determined. This result
generalizes Balinski's theorem on the one-dimensional skeleton of
a $d$-dimensional convex polytope.
\end{abstract}

\maketitle

\section{Introduction}
\label{intro}

The combinatorial theory of convex polytopes has provided
mathematicians with interesting problems since antiquity. Some of
these problems concern skeleta of polytopes and their
connectivity; see for instance \cite{Ka1}. Balinski's theorem
\cite{Ba} (see also Theorem \ref{thm:ba}) gives a sharp lower
bound on the (vertex) connectivity degree of the abstract graph
$\gG(P)$ defined by the one-dimensional skeleton of a
$d$-dimensional convex polytope $P$. This graph can also be
understood as the dual graph of the normal fan $\nN_P$ of $P$ or,
equivalently, as the dual graph of the $(d-1)$-dimensional
skeleton of any polytope $Q$ which is polar-dual to $P$. By
definition, these are the simple graphs on the node set of maximal
faces of $\nN_P$ and of $(d-1)$-dimensional faces (facets) of $Q$,
respectively, in which two such faces are adjacent if they share a
common codimension one face. It is natural to inquire about the
connectivity of the dual graphs of other skeleta of $\nN_P$ or
$Q$. These graphs correspond to the simple graphs $\gG_k (P)$
defined for nonnegative integers $k$ as follows. The nodes of
$\gG_k (P)$ are the $k$-dimensional faces of $P$ and two such
faces are adjacent if there exists a $(k+1)$-dimensional face of
$P$ which contains them both. In other words, $\gG_k (P)$ is the
graph on the node set of rank $k+1$ elements in the face lattice
of $P$, two such elements being adjacent if they have a common
cover in this lattice. For $k=0$, the graph $\gG_k (P)$ reduces to
the graph $\gG(P)$ which appears in Balinski's theorem. It is
folklore that the graphs $\gG_k (P)$ are connected; see
\cite[Theorem 19.5.2]{Ka1}. Their higher (vertex) connectivity is
the subject of this paper.

Let $m$ be a positive integer and recall that an abstract graph
$\gG$ is said to be $m$-connected if $\gG$ has at least $m+1$
nodes and any graph obtained from $\gG$ by deleting $m-1$ or fewer
nodes and their incident edges is connected. Our main result is as
follows.
\begin{theorem}
For fixed nonnegative integers $k$ and $d$ satisfying $0 \le k \le
d-1$, let $m_k (d)$ denote the largest integer $m$ such that the
graph $\gG_k (P)$ is $m$-connected for all convex polytopes $P$ of
dimension $d$. We have
\[ m_k (d) \, = \, \begin{cases}
d, & \text{if \ $k=d-2$} \\
(k+1) (d-k), & \text{otherwise.}
\end{cases} \]
\label{thm0}
\end{theorem}

A few remarks on Theorem \ref{thm0} are in order. Any node of an
$m$-connected graph $\gG$ must have at least $m$ neighbors, since
the graph obtained from $\gG$ by removing all neighbors of a given
node either has a single node or else is disconnected. Theorem
\ref{thm0} is made plausible by the fact that any $k$-dimensional
face $F$ of a $d$-dimensional polytope $P$ has at least $(k+1)
(d-k)$ neighbors in $\gG_k (P)$. Indeed, such a face $F$ is
contained in at least $d-k$ faces of $P$ of dimension $k+1$, each
one of those has at least $k+1$ faces of dimension $k$ other than
$F$ and all these $k$-dimensional faces of $P$ are pairwise
distinct and are neighbors of $F$ in $\gG_k (P)$. On the other
hand, if $P$ is a $d$-dimensional simplex then any $k$-dimensional
face of $P$ has exactly $(k+1) (d-k)$ neighbors in $\gG_k (P)$ and
therefore $\gG_k (P)$ is not $m$-connected for any value of $m$
which exceeds $(k+1) (d-k)$. This example shows that $m_k (d) \le
(k+1) (d-k)$ for all $k$ and $d$. Theorem \ref{thm0} reduces to
Balinski's theorem in the case $k=0$ and is trivial for $k=d-1$,
since $\gG_{d-1} (P)$ is the complete graph on the set of facets
of $P$ and any $d$-dimensional polytope has at least $d+1$ facets.

The connectivity of skeleta of polytopes and more general cell
complexes was studied by Fl$\o$ystad \cite{Fl} from a homological
point of view and, more recently, by Bj\"orner \cite{Bj2} from a
homotopy-theoretic point of view. Fl$\o$ystad showed that
Balinski's theorem holds for one-dimensional skeleta of a class of
cell complexes which includes all finite Cohen-Macaulay polyhedral
complexes, namely that of finite Cohen-Macaulay regular cell
complexes with the intersection property. It is plausible  (see
Conjecture \ref{conj}) that Theorem \ref{thm0} generalizes for $0
\le k \le d-2$ in this respect.

This paper is structured as follows. Section \ref{sec:pre} reviews
graph-theoretic terminology and basic background on convex
polytopes. The special cases $k=d-2$ and $k=1$ of Theorem
\ref{thm0} are proved in Sections \ref{sec:ridgets} and
\ref{sec:edges}, respectively. The proof of the theorem is
completed in Section \ref{sec:proof}. Section \ref{sec:cell}
discusses possible generalizations to classes of cell complexes.
\begin{remark} {\rm
After this paper was written, it was pointed out to the author by
Ronald Wotzlaw that the graphs $\gG_k (P)$ have been previously
studied by Sallee \cite{Sa}. In the notation of Theorem
\ref{thm0}, it is proved in \cite{Sa} (see equation (7.18) on page
495) that $(k+1) (d-k) - k \le m_k (d) \le (k+1) (d-k)$. More
generally, given integers $0 \le r < s \le d-1$, upper and lower
bounds are given in \cite{Sa} for the connectivity degree of the
graph on the node set of $r$-dimensional faces of a
$d$-dimensional convex polytope $P$, in which two such faces are
adjacent if there exists an $s$-dimensional face of $P$ which
contains them both (thus our setting corresponds to the case $r=k$
and $s=k+1$). More precise results are obtained in \cite{Sa} for
various other notions of connectivity for the incidence graphs
between faces of dimension $r$ and faces of dimension $s$ of
$d$-dimensional polytopes. \qed}
\label{rem0}
\end{remark}

\section{Preliminaries}
\label{sec:pre}

All graphs considered in this paper will be simple (without loops
or multiple edges) and finite. Thus every edge of a graph $\gG$
connects two distinct nodes of $\gG$, called its \emph{endpoints},
and is said to be \emph{incident} to each of these nodes. Two
nodes of $\gG$ are said to be \emph{adjacent} if they are
connected by an edge in $\gG$. A \emph{walk} of length $n$ in
$\gG$ is an alternating sequence $w = (v_0, e_1, v_1,\dots,e_n,
v_n)$ of nodes and edges, such that $v_{i-1}$ and $v_i$ are the
endpoints of $e_i$ for $1 \le i \le n$. We say that $w$
\emph{connects} nodes $v_0$ and $v_n$, which are the
\emph{endpoints} of $w$. Thus $\gG$ is connected if any two nodes
can be connected by a walk in $\gG$. Given a subset $V$ of the set
of nodes of $\gG$, we denote by $\gG \sm V$ the graph obtained
from $\gG$ by deleting the nodes in $V$ and all incident to these
nodes edges of $\gG$. Given a positive integer $m$, the graph
$\gG$ is said to be \emph{$m$-connected} if it has at least $m+1$
nodes and $\gG \sm V$ is connected for all subsets $V$ of the set
of nodes of $\gG$ with cardinality at most $m-1$.

A \emph{convex polytope} $P$ is defined as the convex hull of a
finite set of points in $\RR^N$. The \emph{dimension} of $P$ is
the dimension of the affine hull of $P$, as an affine subspace of
$\RR^N$. A \emph{face} of $P$ is either a subset of $P$ on which
some linear functional on $\RR^N$ achieves its minimum on $P$ or
else the empty set. Any face of $P$ is a polytope in $\RR^N$ and
the intersection of two faces of $P$ is again a face of $P$. Faces
of $P$ of dimension zero or one are called \emph{vertices} or
\emph{edges}, respectively, and faces of codimension one are
called \emph{facets}. Every edge has exactly two vertices, also
called its \emph{endpoints}. The following lemma is well-known;
see, for instance \cite[Section 1.2]{Ka2}.
\begin{lemma}
Any $d$-dimensional convex polytope has at least ${d+1 \choose
i+1}$ faces of dimension $i$ for all nonnegative integers $i$.
\label{prop:LBT}
\end{lemma}

The \emph{graph of $P$}, denoted by $\gG(P)$, is the abstract
graph which has as nodes the vertices of $P$ and in which two
nodes are adjacent if they are the endpoints of an edge of $P$. We
will need the following slightly stronger version of Balinski's
theorem, which follows, for instance, from the proof given in
\cite[Section 3.5]{Zi}.

\begin{theorem}
{\rm (Balinski \cite{Ba})} Given a $d$-dimensional convex polytope
$P \subset \RR^N$, the graph $\gG(P) \sm V$ is connected for any
subset $V$ of the vertex set of $P$ which is contained in some
$(d-2)$-dimensional affine subspace of $\RR^N$. In particular,
$\gG (P)$ is $d$-connected.
\label{thm:ba}
\end{theorem}

For any convex polytope $P$ there exists a polytope $Q$
(necessarily of the same dimension) whose set of faces is in
inclusion-reversing bijection with the set of faces of $P$. Such a
polytope is said to be \emph{polar-dual} to $P$. Given an
$r$-dimensional face $F$ of a $d$-dimensional convex polytope $P$,
the set of faces of $P$ which contain $F$ is in
inclusion-preserving bijection with the set of faces of a
$(d-r-1)$-dimensional polytope $P_F$, called a \emph{face figure}
(or \emph{vertex figure} if $F$ is a vertex) of $P$ at $F$. It is
clear that walks in the graph $\gG_k (P_F)$ (as defined in the
introduction) correspond bijectively to walks in $\gG_{k+r+1} (P)$
which involve only faces of $P$ containing $F$. For more
information on convex polytopes and their combinatorial structure
we refer the interested reader to \cite{Ka2, Zi}.

\section{The case $k=d-2$}
\label{sec:ridgets}

In this section we restate and prove Theorem \ref{thm0} in the
case $k=d-2$.

\begin{proposition}
The graph $\gG_{d-2} (P)$ is $d$-connected for all convex
polytopes $P$ of dimension $d$. Moreover, in any dimension $d \ge
2$ there exist convex polytopes for which $\gG_{d-2} (P)$ is not
$(d+1)$-connected.
\label{prop:i=d-2}
\end{proposition}

It is perhaps easier to visualize the graph $\gG_{d-2} (P)$ in
terms of a polytope $Q$ which is polar-dual to $P$. Since faces of
$P$ of dimension $d-2$ and $d-1$ correspond to edges and vertices
of $Q$, respectively, $\gG_{d-2} (P)$ is isomorphic to the graph
$\Gamma(Q)$ on the node set of edges of $Q$ in which two nodes are
adjacent if they have a vertex of $Q$ as a common endpoint.

\medskip
\noindent \emph{Proof of Proposition \ref{prop:i=d-2}.} Let $Q$ be
a polytope polar-dual to $P$, so that $\gG_{d-2} (P)$ may be
replaced by the graph $\Gamma(Q)$ on the node set of edges of $Q$.
To show that $\Gamma(Q)$ is $d$-connected, let $\eE$ be any subset
of the set of edges of $Q$ of cardinality at most $d-1$. We
observe first that given any $e \in \eE$, it is possible to choose
some two-dimensional face $F$ of $P$ containing $e$, so that $e$
is the unique edge of $F$ which belongs to $\eE$. Indeed, in view
of our assumption on the cardinality of $\eE$, this is so because
$e$ is contained in at least $d-1$ faces of $P$ of dimension 2 and
any two such faces share no edge other than $e$ in common. Since
$\gG(Q)$ is connected, so is $\Gamma(Q)$. Our previous remark
shows that any walk $w$ in $\Gamma(Q)$ connecting two nodes not in
$\eE$ can be transformed to one connecting the same nodes that
does not involve elements of $\eE$. This can be done by replacing
any node $e \in \eE$ that may appear in $w$ with the sequence of
edges other than $e$ (and vertices other than the endpoints of
$e$), ordered appropriately, of a two-dimensional face of $P$
which contains $e$ but no other edge in $\eE$. It follows that
$\Gamma(Q) \sm \eE$ is connected and hence that $\Gamma(Q)$ is
$d$-connected.

To prove the second statement in the proposition let $I$ be a line
segment, let $Q = \Delta \times I$ be the prism over a
$(d-1)$-dimensional simplex $\Delta$ and denote by $\eE$ the set
of edges of $Q$ of the form $v \times I$, where $v$ is a vertex of
$\Delta$. It is clear that the set $\eE$ has $d$ elements and that
the graph $\Gamma(Q) \sm \eE$ is disconnected. This implies that
$\Gamma(Q)$ is not $(d+1)$-connected. \qed

\section{The case $k=1$}
\label{sec:edges}

In this section we prove Theorem \ref{thm0} in the case $k=1$ as
follows.

\begin{proposition}
The graph $\gG_1 (P)$ is $(2d-2)$-connected for all convex
polytopes $P$ of dimension $d \ge 4$.
\label{thm1}
\end{proposition}
\begin{proof}
Recall that the set of nodes of the graph $\gG_1 (P)$ coincides
with the edge set of $P$. Let $\eE$ be any subset of this set of
cardinality less than $2d-2$ and let $f$ and $g$ be any two edges
of $P$ not in $\eE$. We need to show that $f$ and $g$ can be
connected by a walk in the graph $\gG_1 (P) \sm \eE$. For any
vertex $v$ and any edge $e$ of $P$, we will denote by $s(v)$ the
number of edges in $\eE$ which have $v$ as an endpoint and by
$t(e)$ the number of edges in $\eE$ which have at least one common
endpoint with $e$.

\medskip
\noindent {\bf Claim:} There exists a walk $w = (v_0, e_1,
v_1,\dots,e_n, v_n)$ in $\gG (P)$ such that $v_0$ is an endpoint
of $f$, $v_n$ is an endpoint of $g$ and the following hold:
\begin{itemize}
\itemsep=0pt
\item[{\rm (a)}] $s(v_i) \le d-2$ for $0 \le i \le n$
and

\item[{\rm (b)}] $t(e_i) \le d-1$ for all $1 \le i \le n$ with
$e_i \in \eE$.
\end{itemize}

\noindent Given the claim, we can proceed as follows. For any
index $1 \le i \le n$ with $e_i \in \eE$, there exist at least
$d-1$ two-dimensional faces of $P$ which contain $e_i$ and any two
of them have no edge other than $e_i$ in common. Because of
condition (b), in at least one of these faces none of the two
edges which share exactly one common endpoint with $e_i$ belongs
to $\eE$. Therefore we may choose edges $g_i$ and $f_i$ of $P$ not
in $\eE$ which are adjacent nodes in the graph $\gG_1 (P)$, so
that $v_{i-1}$ is an endpoint of $g_i$ and $v_i$ is an endpoint of
$f_i$. We set $g_i = f_i = e_i$ for those indices $1 \le i \le n$
for which $e_i$ is not an element of $\eE$. We also set $f_0 = f$
and $g_{n+1} = g$. For $0 \le i \le n$ we observe that $f_i$ and
$g_{i+1}$ are edges of $P$ not in $\eE$ which have $v_i$ as a
common endpoint and denote by $P_i$ the vertex figure of $P$ at
$v_i$ and by $V_i$ the set of vertices of $P_i$ which correspond
to elements of $\eE$ having $v_i$ as an endpoint. Since, by
Theorem \ref{thm:ba}, the graph $\gG(P_i)$ is $(d-1)$-connected
and, by condition (a), the set $V_i$ has no more than $d-2$
elements, the vertices of $P_i$ which correspond to $f_i$ and
$g_{i+1}$ can be connected by a walk in $\gG(P_i) \sm V_i$. This
implies that $f_i$ and $g_{i+1}$ can be connected by a walk in the
graph $\gG_1 (P) \sm \eE$ for $0 \le i \le n$. Therefore $f_0=f$
and $g_{n+1}=g$ can also be connected by a walk in this graph.

Thus it suffices to prove the claim. Let us call a vertex $v$ of
$P$ bad if it is an endpoint of at least $(d+1)/2$ edges in $\eE$
and good otherwise. We will also call an edge $e$ of $P$ bad if it
violates condition (b), meaning that $e \in \eE$ and $t(e) \ge d$.
Otherwise we call $e$ good. Clearly, any vertex of $P$ which
violates condition (a) is bad and any bad edge has at least one
bad endpoint. As a result, any walk in $\gG(P)$ which does not go
through bad vertices satisfies conditions (a) and (b). Let $q$
denote the cardinality of $\eE$. Since each edge of $P$ has only
two endpoints, the number $p$ of bad vertices of $P$ satisfies $p
\left\lceil (d+1)/2 \right\rceil \le 2q \le 2(2d-3)$. From this
and our assumption $d \ge 4$ it follows that $p \le d-1$, that is
there exist at most $d-1$ bad vertices. Therefore, by Theorem
\ref{thm:ba}, deleting all bad vertices of $P$ from $\gG(P)$ and
their incident edges results in a connected graph. Thus the claim
will follow if we can show that each of the edges $f$ and $g$
either has a good endpoint or else one of its endpoints, say $v$,
satisfies $s(v) \le d-2$ and is connected to a good vertex by a
good edge of $P$. We will prove this statement for $f$, the same
arguments applying for $g$. Let $a$ and $b$ be the endpoints of
$f$. We distinguish two cases:

{\bf Case 1:} At least one of the endpoints of $f$, say $a$,
satisfies $s(a) \ge d-1$. Since $\eE$ has at most $2d-3$ elements
and $f$ is not in $\eE$, we must have $s(b) \le d-2$. As a result,
there exists at least one edge $e$ of $P$ other than $f$ which has
$b$ as an endpoint and does not belong to $\eE$. Simple counting
shows that at least one of the endpoints of $e$ is good. Since $e$
is a good edge, we are done in this case.

{\bf Case 2:} We have $s(a) \le d-2$ and $s(b) \le d-2$. As a
result, each of $a, b$ is an endpoint of at least one edge of $P$
not in $\eE$, other than $f$. There is nothing to prove if at
least one of $a, b$ is good, so we may assume that both of them
are bad. Suppose first that there exist distinct vertices $a'$ and
$b'$ of $P$ other than $a, b$ which are connected to $a$ and $b$,
respectively, by edges of $P$ not in $\eE$. Once again, simple
counting shows that at least one of $a', b'$ must be good and the
desired statement follows. Otherwise $a$ and $b$ must be connected
to a vertex $c$ of $P$ with edges not in $\eE$ and we must have
$s(a) = s(b) = d-2$. It follows that $s(c) \le 1$ and that $c$ is
good, so we are done in this case as well.
\end{proof}

\section{Proof of Theorem \ref{thm0}}
\label{sec:proof}

In this section we complete the proof of our main theorem. The
structure of the proof is similar to that of Proposition
\ref{thm1}.

\medskip
\noindent \emph{Proof of Theorem \ref{thm0}.} Let us write $n_k
(d) = (k+1)(d-k)$. We have already remarked in the introduction
that there exist $d$-dimensional convex polytopes $P$ such that
$\gG_k (P)$ is not $m$-connected for $m > n_k (d)$ and that
Theorem \ref{thm0} is trivial for $k=d-1$. In view of Proposition
\ref{prop:i=d-2}, it remains to show that for $1 \le k \le d-3$
the graph $\gG_k (P)$ is $n_k (d)$-connected for all convex
polytopes $P$ of dimension $d$.

We proceed by induction on $d$ and $k$, where the case $k=1$ was
treated by Proposition \ref{thm1}. Assume that $k \ge 2$. Let $U$
be any subset of the set of $k$-dimensional faces of $P$ of
cardinality less than $n_k (d)$ and let $F$ and $G$ be two
$k$-dimensional faces of $P$ not in $U$. We need to show that $F$
and $G$ can be connected by a walk in $\gG_k (P) \sm U$.

\medskip
\noindent {\bf Claim:} There exists a walk $w$ in $\gG (P)$ which
connects a vertex of $F$ to a vertex of $G$ and has the following
properties:
\begin{itemize}
\itemsep=0pt \item[{\rm (a)}] no edge of $w$ is contained in ${d-1
\choose k-1}$ or more faces in $U$ and

\item[{\rm (b)}] no node of $w$ belongs to $n_{k-1}(d-1)$ or more
faces in $U$.
\end{itemize}

\noindent Given the claim, the proof proceeds as follows. Let $w =
(v_0, e_1, v_1,\dots,e_n, v_n)$ be a walk as in the claim and set
$F_0=F$ and $F_{n+1} = G$. It follows from Lemma \ref{prop:LBT}
that each edge of $P$ is contained in at least ${d-1 \choose k-1}$
faces of $P$ of dimension $k$. Therefore, in view of our condition
(a), for each index $1 \le i \le n$ we may choose a
$k$-dimensional face $F_i$ of $P$ not in $U$ which contains the
edge $e_i$. Note that $v_i$ is a vertex of both $F_i$ and
$F_{i+1}$ for $0 \le i \le n$. Let $P_i$ denote the vertex figure
of $P$ at $v_i$ and let $U_i$ denote the set of
$(k-1)$-dimensional faces of $P_i$ which correspond to the faces
of $U$ containing $v_i$. By the induction hypothesis on $k$, the
graph $\gG_{k-1} (P_i)$ is $n_{k-1}(d-1)$-connected. By condition
(b), this implies that $\gG_{k-1} (P_i) \sm U_i$ is connected and
hence that $F_i$ and $F_{i+1}$ can be connected by a walk in the
graph $\gG_k (P) \sm U$. Since this holds for all $0 \le i \le n$,
we conclude that $F_0=F$ and $F_{n+1}=G$ can be connected by a
walk in $\gG_k (P) \sm U$ as well. It follows that $\gG_k (P) \sm
U$ is connected and hence that $\gG_k (P)$ is $n_k (d)$-connected,
as desired. It therefore suffices to prove the claim. We
distinguish two cases:

{\bf Case 1:} $k=2$. We are given that $d \ge 5$ and that $U$ is a
set of two-dimensional faces of $P$ of cardinality less than $n_k
(d) = 3d-6$ and note that $n_{k-1} (d-1) = 2d-4$. Let us call an
edge or vertex of $P$ bad if this edge or vertex is contained in
at least $d-1$ or $2d-4$, respectively, elements of $U$. The
following hold: (i) there exist at most two bad edges of $P$, (ii)
there exist at most two bad vertices of $P$ and (iii) if $v$ is a
bad vertex of $P$ and $e$ is a bad edge, then $v$ is an endpoint
of $e$. Indeed, the existence of three bad edges of $P$ would
require at least $3d-6$ elements of $U$, since given any two edges
of a polytope $P$, there exists at most one 2-dimensional face of
$P$ which contains both of these edges. In view of our assumption
on the cardinality of $U$, this proves (i). We next observe that
if $u$ and $v$ are distinct bad vertices of $P$, then there exist
at least $d-1$ elements of $U$ which contain both $u$ and $v$.
Therefore any two bad vertices are connected by a bad edge and
(ii) follows from (i). Finally, if $v$ is a bad vertex of $P$ and
$e$ is a bad edge, then there must be at least two elements of $U$
which contain both $v$ and $e$. The intersection of these has to
equal $e$ and contains $v$. This proves (iii). It follows from
facts (i)--(iii) that one can choose a set $V$ consisting of at
most two vertices of $P$ such that $\gG (P) \sm V$ contains no bad
vertex or edge. This completes the proof of the claim in this case
since $\gG (P) \sm V$ is connected, by Theorem \ref{thm:ba}, and
any walk in this graph connecting a vertex of $F$ to a vertex of
$G$ satisfies conditions (a) and (b).

{\bf Case 2:} $k \ge 3$. Since we have ${d-1 \choose k-1} \ge k
(d-k) = n_{k-1}(d-1)$ for $d \ge k+3 \ge 6$, condition (a) follows
from (b) and can thus be ignored. Let $V$ be the set of vertices
of $P$ which belong to at least $n_{k-1}(d-1) = k(d-k)$ faces in
$U$. We will show that any $k+1$ vertices, say $v_0,
v_1,\dots,v_k$, in $V$ are affinely dependent. Indeed, since $U$
has less than $n_k (d) = (k+1)(d-k)$ elements and each $v_i$
belongs to at least $k(d-k)$ of them, there must be at least $k$
elements of $U$ which contain all of $v_0, v_1,\dots,v_k$.
Clearly, the intersection of any two if these $k$ elements
contains the $v_i$ and has affine dimension at most $k-1$, so
$v_0, v_1,\dots,v_k$ must be affinely dependent. It follows that
the dimension of the affine span of $V$ is at most $k-1$. As a
consequence, each one of the $k$-dimensional faces $F$ and $G$ has
at least one vertex not in $V$ and, by Theorem \ref{thm:ba}, the
graph $\gG (P) \sm V$ is connected. These two facts imply the
existence of a walk $w$ in $\gG (P)$ with the claimed properties.
\qed

\section{Cell complexes}
\label{sec:cell}

In this section we discuss possible generalizations of Theorem
\ref{thm0} to certain classes of regular cell complexes. We will
assume some familiarity with regular cell complexes and standard
notions in topological combinatorics; excellent sources on these
topics are \cite[Section 4.7]{OM} and the article \cite{Bj}.

A \emph{regular cell complex} is a finite collection $\cC$ of
balls in a Hausdorff space $\|\cC\| = \bigcup_{\sigma \in \cC} \,
\sigma$, called \emph{cells} or \emph{faces}, such that:
\begin{itemize}
\itemsep=0pt
\item[{\rm (i)}] $\varnothing \in \cC$,

\item[{\rm (ii)}] the relative interiors of the nonempty cells
partition $\|\cC\|$ and

\item[{\rm (iii)}] the boundary of any cell in $\cC$ is a union of
cells in $\cC$.
\end{itemize}

\noindent Cells of dimension zero or one are called
\emph{vertices} or \emph{edges}, respectively, and cells which are
maximal with respect to inclusion are called \emph{facets}. The
dimension of $\cC$ is the maximum dimension of a cell. The
(loop-free) abstract graph $\gG(\cC)$ defined by the vertices and
edges of $\cC$ is called the \emph{graph} of $\cC$.

A \emph{polyhedral complex} in $\RR^N$ is a regular cell complex
each of whose cells is a polytope in $\RR^N$. A \emph{simplicial
complex} is a polyhedral complex in which every cell is a simplex.
A regular cell complex $\cC$ is said to have the
\emph{intersection property} if the intersection of any two cells
in $\cC$ is also a cell in $\cC$ (in particular, the graph
$\gG(\cC)$ is simple). For instance, polyhedral complexes have the
intersection property. A regular cell complex $\cC$ with the
intersection property is said to be \emph{Cohen-Macaulay} (over a
field $\KK$) \cite[Section 2]{Fl} if, under the inclusion partial
order, it is a Cohen-Macaulay poset (over $\KK$); see
\cite[Section 11]{Bj} for the notion of Cohen-Macauliness for
simplicial complexes and posets. Such a complex $\cC$ is pure,
meaning that all facets of $\cC$ have the same dimension, and
strongly connected, meaning that for any two facets $\tau$ and
$\tau'$ of $\cC$ there exists a sequence $\tau=\tau_0,
\tau_1,\dots,\tau_n = \tau'$ of facets of $\cC$ such that
$\tau_{i-1}$ and $\tau_i$ intersect on a common face of
codimension one for all $1 \le i \le n$. The following
generalization of Balinski's theorem was proved by Fl$\o$ystad, as
a consequence of \cite[Corollary 2.7]{Fl}.

\begin{theorem}
{\rm (Fl$\o$ystad \cite{Fl})} For any $d$-dimensional
Cohen-Macaulay regular cell complex with the intersection
property, the graph $\gG(\cC)$ is $d$-connected.
\label{thm:fl}
\end{theorem}
Let $\gG_k (\cC)$ denote the simple graph on the node set of
$k$-dimensional cells of $\cC$ in which two such cells are
adjacent if there exists a $(k+1)$-dimensional cell of $\cC$ which
contains them both. Theorem \ref{thm:fl} is the special case $k=0$
of the following statement.
\begin{conjecture}
For any $d$-dimensional Cohen-Macaulay regular cell complex $\cC$
with the intersection property, the graph $\gG_k (\cC)$ is
\begin{itemize}
\itemsep=0pt
\item[$\circ$] $(k+1)(d-k)$-connected if $0 \le k \le
d-3$ and

\item[$\circ$] $d$-connected if $k=d-2$.
\end{itemize}
\label{conj}
\end{conjecture}
The $d$-dimensional simplicial complex which has two
$d$-dimensional simplices as facets, intersecting on a common
codimension one face, shows that Theorem \ref{thm0} does not
extend to the setup of Conjecture \ref{conj} for $k=d-1$. We can
verify this conjecture in the special cases which appear in the
following statement.
\begin{proposition}
Under the assumptions of Conjecture \ref{conj}, the graph $\gG_k
(\cC)$ is
\begin{itemize}
\itemsep=0pt
\item[(i)] $(k+1)(d-k)$-connected if $0 \le k \le
d-3$ and $\cC$ is a polyhedral complex,

\item[(ii)] $d$-connected if $k=d-2$.
\end{itemize}
\label{thm:cell}
\end{proposition}
\begin{proof}
Suppose first that $\cC$ is a polyhedral complex and that $0 \le k
\le d-3$ (a similar argument works for $k=d-2$). Let $U$ be any
subset of the set of $k$-dimensional faces of $\cC$ of cardinality
less than $(k+1)(d-k)$ and let $F$ and $G$ be two $k$-dimensional
faces of $\cC$ not in $U$. We need to show that these two faces
can be connected by a walk in the graph $\gG_k (\cC) \sm U$. Let
$P$ and $Q$ be facets of $\cC$ (necessarily of dimension $d$)
containing $F$ and $G$, respectively. By strong connectivity of
$\cC$, we may choose a sequence $P=P_0, P_1,\dots,P_n = Q$ of
facets of $\cC$ such that $P_{i-1}$ and $P_i$ intersect on a
common $(d-1)$-dimensional face for all $1 \le i \le n$. By Lemma
\ref{prop:LBT}, each intersection $P_{i-1} \cap P_i$ has at least
${d \choose k+1}$ faces of dimension $k$. Since ${d \choose k+1}
\ge (k+1)(d-k)$ for $k \le d-3$, there exists at least one
$k$-dimensional face, say $F_i$, of $P_{i-1} \cap P_i$ which is
not an element of $U$. We let $F_0 = F$ and $F_{n+1} = G$ and note
that $F_i$ and $F_{i-1}$ can be connected with a walk in $\gG_k
(P_{i-1}) \sm U$ for all $1 \le i \le n+1$, by Theorem \ref{thm0}.
It follows that $F$ and $G$ can be connected with a walk in $\gG_k
(\cC) \sm U$. This proves (i).

Suppose now that $k=d-2$ and that $\cC$ is as in Conjecture
\ref{conj}. We can proceed as in the proof of Proposition
\ref{prop:i=d-2} with no need to pass to a dual to $\cC$ object.
Indeed, let $\eE$ be any subset of the set of $(d-2)$-dimensional
faces of $\cC$ of cardinality at most $d-1$. We observe first that
any $e \in \eE$ has at least $d-1$ codimension one faces and that
no $(d-2)$-dimensional face of $\cC$ other than $e$ contains two
or more of those. Therefore, it is possible to choose a
codimension one face $\sigma$ of $e$ so that $e$ is the unique
element of $\eE$ which contains $\sigma$. We then check that in
any part of a walk in $\gG_{d-2} (\cC)$ of the form $(\tau, e,
\tau')$, where $e \in \eE$ and $\tau, \tau'$ are faces of
dimension $d-1$ containing $e$, the node $e$ can be replaced by a
walk that does not involve elements of $\eE$ as follows. Since the
inclusion poset of faces of $\cC$ which strictly contain $e$ is
Cohen-Macaulay of rank one, and hence connected, we may assume
that some facet $\rho$ of $\cC$ contains both $\tau$ and $\tau'$.
We pick a codimension one face $\sigma$ of $e$ as above and
observe that the set of faces of $\cC$ containing $\sigma$ and
contained in $\rho$ is in inclusion preserving bijection with the
set of faces of a polygon $\Pi$ (this holds more generally for
Gorenstein* posets of rank 3) and that $e$ and $\tau, \tau'$
correspond to a vertex and its two incident edges in $\Pi$. Hence
we can deviate our given walk in $\gG_{d-2} (\cC)$ around $e$
through the boundary of $\Pi$, thus avoiding nodes in $\eE$ by the
defining property of $\sigma$. To complete the proof of (ii) it
remains to comment that the graph $\gG_{d-2} (\cC)$ is connected.
This holds because the inclusion order on the set of faces of
$\cC$ of dimension $d-2$ or $d-1$ inherits the Cohen-Macaulay
property from $\cC$.
\end{proof}

\medskip
\noindent \emph{Acknowledgements}: The author thanks Bernd
Sturmfels for encouraging discussions, Vic Reiner for providing an
example of a regular cell complex with a nonregular face figure
and Ronald Wotzlaw for the content of Remark \ref{rem0} and for
useful pointers to the literature.


\begin{thebibliography}{99}
%
\bibitem{Ba}
M.L.~Balinski,
\emph{On the graph structure of convex polyhedra in
$n$-space},
Pacific J. Math. {\bf~11} (1961), 431--434.
%
\bibitem{Bj}
A.~Bj\"orner,
\emph{Topological methods},
in \emph{Handbook of Combinatorics} (R.L.~Graham, M.~Gr\"otschel
and L.~Lov\'asz, eds.), North Holland, Amsterdam, 1995,
pp.~1819--1872.
%
\bibitem{Bj2}
A.~Bj\"orner,
\emph{Connectivity of polytopes and Cohen-Macaulay rigidity},
in preparation.
%
\bibitem{OM}
A.~Bj\"orner, M.~Las~Vergnas, B.~Sturmfels, N.~White and
G.M.~Ziegler,
Oriented Matroids,
Encyclopedia of Mathematics and Its Applications {\bf~46}, Cambridge
University Press, Cambridge, 1993; second printing, 1999.
%
\bibitem{Fl}
G.~Fl$\o$ystad,
\emph{Cohen-Macaulay cell complexes},
in \emph{Algebraic and Geometric Combinatorics} (C.A.~Athanasiadis,
V.V.~Batyrev, D.I.~Dais, M.~Henk and F.~Santos, eds.), Proceedings
of a Euroconference in Mathematics, August 2005, Anogia, Crete,
Greece, Contemporary Mathematics {\bf~423}, American Mathematical
Society, Providence, RI, 2007, pp.~205--220.
%
\bibitem{Ka1}
G.~Kalai, \emph{Polytope skeletons and paths},
in \emph{Handbook of Discrete and Computational Geometry} (J.E.~Goodman
and J.~O'Rourke, eds.), CRC Press, Boca-Raton, FL, 1997, pp.~331--344.
%
\bibitem{Ka2}
G.~Kalai,
\emph{Some aspects of the combinatorial theory of convex
polytopes},
in \emph{Polytopes: abstract, convex and computational} (Scarborough, ON,
1993), Kluwer Academic Publishers, Dordrecht, 1994, pp.~205--229.
%
\bibitem{Sa}
G.T.~Sallee,
\emph{Incidence graphs of convex polytopes},
J.
Combin. Theory {\bf~2} (1967), 466--506.
%
\bibitem{Zi}
G.M.~Ziegler,
Lectures on Polytopes,
Graduate Texts in Mathematics {\bf~152}, Springer-Verlag, New York, 1995.
%
\end{thebibliography}
\end{document}